\documentclass{article}

\usepackage{amsmath}
\usepackage{amssymb}
\usepackage[margin=1in]{geometry}
\usepackage{graphicx}
\usepackage{hyperref}
\usepackage{tikz}
\usepackage{caption}

\newcommand{\R}{\mathbb{R}}

\newcommand{\dd}{\, \mathrm{d}}
\newcommand{\LL}{\mathcal{L}}

\usepackage{amsthm}

\theoremstyle{definition}
\newtheorem{theorem}{Theorem}
\newtheorem{lemma}[theorem]{Lemma}

\theoremstyle{definition}

\newtheorem{remark}[theorem]{Remark}

\title{The entropy production is not always monotone in the space-homogeneous Boltzmann equation}
\author{Luis Silvestre}
\date{\today}

\begin{document}

\maketitle

\begin{abstract}
    We show an example of a function and a collision kernel for which the entropy production increases in time when we flow it by the space-homogeneous Boltzmann equation. The collision kernel is not any of the physically motivated kernels that are commonly used in the literature. In this particular setting, our result disproves a conjecture of McKean from 1966.
\end{abstract}  

\section{Introduction}
If a function $f = f(t,v)$ solves the space-homogeneous Boltzmann equation, then the entropy is monotone decreasing in time. Here, we use the following expression for the entropy (sometimes the opposite sign convention is used):
    \[ H(f) = -\int_{\R^d} f \log f \, dv. \]
The entropy production is defined as
    \[ D(f) = \frac{d}{dt} H(f) = -\int_{\R^d} Q(f,f) \log f \, dv. \]
In 1966, McKean conjectured \cite{mckean1966speed} that the entropy production is also monotone decreasing in time. This conjecture has been tested numerically \cite{filbet2006solving, filbet2003high}. While oscillations in the entropy production are possible in the space-inhomogeneous setting, they have never been observed in the space-homogeneous Boltzmann equation. This issue is mentioned in \cite{desvillettes2005trend,filbet2006solving}. The problem is also mentioned in the recent review \cite[Bibliographical notes of Section 3]{villani2025fisher} and in \cite[Chapter 4, Section 4.3]{villani2002review}.

A similar question is whether the Fisher information is monotone decreasing in time. It has recently been verified for the Landau equation \cite{guillen2025landau} as well as for the Boltzmann equation \cite{imbert2026monotonicity}. In both cases, the monotonicity of the Fisher information is proved under some \emph{reasonable} extra assumptions on the operator, which are satisfied in all physically relevant scenarios.

In this paper, we construct an example of a function $f$ and a collision kernel for which the entropy production increases in time when we flow it by the space-homogeneous Boltzmann equation. The collision kernel is rather special. It is not any of the physically motivated ones. It does not satisfy the assumptions for the monotonicity of the Fisher information recently developed in \cite{imbert2026monotonicity}. The result disproves McKean's conjecture in this unusual setting.

We recall the definition of the Boltzmann collision operator. For a function $f = f(v)$, we have
    \begin{equation} \label{e:Q}
        Q(f,f)(v) = \int_{\R^d} \int_{\mathbb{S}^{d-1}} (f' f'_* - f f_*) B(|v-v_*|, \cos \theta) \, d\sigma dv_*.
    \end{equation}
where $f' = f(v')$, $f'_* = f(v'_*)$, and
    \[ v' = \frac{v+v_*}{2} + \frac{|v-v_*|}{2} \sigma, \quad v'_* = \frac{v+v_*}{2} - \frac{|v-v_*|}{2} \sigma. \]

We consider a collision kernel of the form
    \begin{equation} \label{e:kernel}
        B(|v-v_*|, \cos \theta) = \Phi(|v-v_*|) b(\cos \theta),
    \end{equation} 
where $\Phi$ is a function of the relative velocity and $b$ is a function of the angle. Some collision kernels are derived from physical models involving interactions between particles. These are the hard-spheres model, and the inverse power-law models. A variety of other collision kernels are considered in the mathematical literature, where the functions $\Phi$ and $b$ are allowed to be quite general. The collision kernel that we use to prove our result is none of the usual kernels that are derived from physical models.  Whether the monotonicity of the entropy production holds in the space-homogeneous Boltzmann equation with physically relevant collision kernels remains an interesting open question. It is conceivable that some form of McKean's conjecture may still hold with additional assumptions on the collision kernel, similar to the setting in \cite{imbert2026monotonicity}.

We state our main theorem in the two-dimensional setting.

\begin{theorem} \label{t:main}
    There exists a nonnegative compactly-supported function $f: \R^2 \to [0,\infty)$ and a collision kernel of the form \eqref{e:kernel} such that the entropy production $D(f(t))$ increases in time when we flow $f$ by the space-homogeneous Boltzmann equation with the given collision kernel.
\end{theorem}

The kernel that we use for our computations is quite singular. We take $b(\cos \theta)$ to be a Dirac mass concentrated on $\theta = \pm \pi/2$, and $\Phi(r)$ to be a Dirac mass concentrated at $r = \sqrt 2$. A posteriori, there is no difficulty in approximating this kernel by a smooth kernel and obtain the same result with a non-singular kernel. However, these approximate kernels would still have $\Phi'(r)/\Phi(r)$ very large near $r = \sqrt 2$, departing very much from the assumptions for the monotonicity of the Fisher information in \cite{imbert2026monotonicity}, and from any of the physically relevant kernels commonly used in the literature. In the rest of the paper, we describe the proof of Theorem \ref{t:main} using the singular kernel described above. Even though the singularity of the kernel is not essential for the result, it helps us organize the computations.

It is possible to construct an example like in Theorem \ref{t:main} in higher dimension using similar ideas. The geometry of the construction would be more complicated, so we choose to present it in the two-dimensional setting for simplicity. The localization of the kernel near specific values of the relative velocity and the deviation angle seems to be an important ingredient for the construction to work out.

\subsection{Other bibliographical notes}

In \cite{mckean1966speed}, McKean presented several ambitious conjectures and speculation. His first conjecture in \cite[Section 13]{mckean1966speed} is that $H''(t) \leq 0$, which is what we disprove in this paper. He then speculated that $H^{(n)}(t)$ has a definite sign for all $n \geq 1$. This more ambitious conjecture (often called the ``super H-theorem'') was disproved in the 1980s for sufficiently high order derivatives. See in \cite{ziff1981approach,lieb1982comment,olaussen1982extension,vigfusson1983higher} and references therein. The result of our present work shows that the second derivative of the entropy can change sign, disproving McKean's conjectures at the lowest possible order. 

Our result also contrasts with the apparent monotonicity of the entropy production observed in numerical tests \cite{pritchard1974entropy,filbet2003high,filbet2006solving}. Naturally, numerical tests are performed with physically relevant collision kernels, and the monotonicity of the entropy production is not disproved in that setting. The result of this paper shows that the monotonicity of the entropy production is not a universal property of the space-homogeneous Boltzmann equation with simply any kernel. It is still conceivable that it may hold under further \emph{reasonable} assumptions.

Interestingly, there is a positive result in a recent preprint by C\^ome Tabary \cite{tabary2026monotonicity}. He shows that the entropy production is monotone decreasing in time for the space-homogeneous Landau equation, in the Maxwell-molecules case, as soon as the directional temperatures are sufficiently  evenly distributed. This is always the case for large enough time. The Maxwell-molecules case corresponds to a collision kernel of the form \eqref{e:kernel} where $\Phi \equiv 1$. It is conceivable that the entropy production may be monotone decreasing in time for the space-homogeneous Boltzmann equation in the case of Maxwell-molecules, even starting from the initial time. The construction in this paper departs sharply from this scenario since we need $\Phi$ to be a Dirac mass, or nearly a Dirac mass.

Understanding the monotonicity of derivatives of the entropy is difficult even for the heat equation, where many open questions remain. See \cite{ledoux2022differentials,wang2025higher}.

\section{The derivative of the entropy production}

In this section, we compute an expression for $\partial_t D(f)$ when $f$ evolves by the space-homogeneous Boltzmann equation. We will later use this expression to find a function $f$ for which $\partial_t D(f) < 0$.

Let us start with the definition of the entropy production:
    \[ D(f) = -\int_{\R^d} Q(f,f) \log f \, dv. \]
It is a simple and standard computation to express $D(f)$ in the following form:
    \[ D(f) = \frac{1}{4} \int_{\R^d} \int_{\R^d} \int_{\mathbb{S}^{d-1}} (f' f'_* - f f_*) \log \left( \frac{ f' \ f'_*}{f \ f_*} \right) B(|v-v_*|, \cos \theta) \, d\sigma dv_* dv. \]
This expression is symmetric in $v$ and $v_*$, and also in $v'$ and $v'_*$. It is also nonnegative, which is the content of the H-theorem.

We proceed to differentiate $D(f)$ in time. We will use the notation $\partial_t f = Q$, which is the space-homogeneous Boltzmann equation. Using the symmetries of the expression for $D(f)$, we obtain
    \begin{align*}
    \partial_t D(f) &= \int_{\R^d} \int_{\R^d} \int_{\mathbb{S}^{d-1}} \left( -Q f_\ast \log \left( \frac{ f' \ f'_*}{f \ f_*} \right) - (f' f'_* - f f_*) \frac Q f \right) B(|v-v_*|, \cos \theta) \, \dd \sigma \dd v_* \dd v \\
    &= - \int_{\R^d} \frac Q f \left( \int_{\R^d} \int_{\mathbb{S}^{d-1}} (f' f'_* - f f_*) B \, \dd \sigma \dd v_* \right) \dd v \\
    &\qquad - \int_{\R^d} \int_{\R^d} \int_{\mathbb{S}^{d-1}} Q f_\ast \log \left( \frac{ f' \ f'_*}{f \ f_*} \right)  B \, \dd \sigma \dd v_* \dd v \\
    &= -\int_{\R^d} \frac{Q^2}{f(v)} \dd v + \int_{\R^d} \int_{\R^d} \int_{\mathbb{S}^{d-1}} Q f_\ast \log \left( \frac{f \ f_*} { f' \ f'_*} \right)  B \, \dd \sigma \dd v_* \dd v
    \end{align*}

We summarize the computation above in the following lemma.

\begin{lemma} \label{l:derivative}
    If $f$ evolves by the space-homogeneous Boltzmann equation, then
    \[ \partial_t D(f) = -\int_{\R^d} \frac{Q^2}{f} \dd v + \int_{\R^d} Q \LL \dd v. \]
    Here $\LL$ is given by
    \[ \LL = \int_{\R^d} \int_{\mathbb{S}^{d-1}} f_\ast \log \left( \frac{f \ f_*} { f' \ f'_*} \right)  B(|v-v_*|, \cos \theta) \, \dd \sigma \dd v_* \dd v. \]
\end{lemma}

The first term in the expression for $\partial_t D(f)$ in Lemma \ref{l:derivative} is negative. The second term does not have a definite sign. 

Testing the expression of $\partial_t D(f)$ in Lemma \ref{l:derivative} with some sample functions, it becomes quite apparent that it is difficult to find a function $f$ for which $\partial_t D(f) > 0$. The first term is negative and it is not clear how to make the second term positive enough to overcome the first term. In the next section, we will find such a function $f$ for a particular collision kernel.

\section{Proof of Theorem \ref{t:main}}

\subsection{The collision kernel}

As we mentioned in the introduction, we will consider a collision kernel of the form \eqref{e:kernel} where $b$ is a Dirac mass concentrated on $\theta = \pm \pi/2$, and $\Phi(r)$ is a Dirac mass concentrated at $r = \sqrt 2$. This kernel $B(|v-v_*|, \cos \theta)$ is a singular measure restricting the integral in the expression \eqref{e:Q} to points where $v$, $v'$, $v_\ast$, $v'_\ast$, in that order, form a square of side length one.

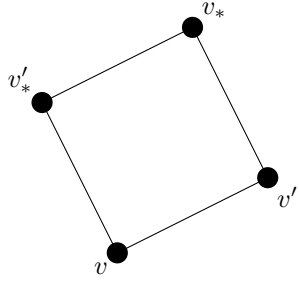
\begin{figure}[ht]
    \centering
    \begin{tikzpicture}[scale=2]
        \draw (0.5,0) -- (1.5,0.5) -- (1,1.5) -- (0,1) -- cycle;
        
        \fill (0.5,0) circle (2pt);
        \fill (1.5,0.5) circle (2pt);
        \fill (1,1.5) circle (2pt);
        \fill (0,1) circle (2pt);
        
        \node[below left] at (0.5,0) {$v$};
        \node[below right] at (1.5,0.5) {$v'$};
        \node[above right] at (1,1.5) {$v_*$};
        \node[above left] at (0,1) {$v'_*$};
        
    \end{tikzpicture}
    \caption{Configuration where vertices form a square of side length one.} \label{f:square}
\end{figure}

Our special choice for $B$ makes the expression for the Boltzmann collision operator $Q$ particularly simple. The integral in \eqref{e:Q} reduces to the following integral on a circle.

\[ Q(f,f)(v) = \int_{S_1} (f' f'_* - f f_*) d\sigma, \]
where (see Figure \ref{f:square})
\begin{align*}
    v' &= v + \sigma, \\
    v'_* &= v + \sigma^\perp, \\
    v_* &= v + (\sigma + \sigma^\perp).
\end{align*}

\subsection{The function}

We proceed to describe the function $f$ for which the entropy production increases in time. To prove Theorem \ref{t:main}, it suffices to find a function $f$ such that $\partial_t D(f) > 0$ if we evolve $f$ by the space-homogeneous Boltzmann equation.

Let $a$ be a positive parameter that will eventually be chosen to be sufficiently
large. Let $\rho>0$ be a small but fixed number, for example we can take $\rho = 1/10$. Let $c>0$ be a small fixed constant to be determined below. We define $f$ as follows.
\begin{align*}
    f(v) &= c a^2 \quad \text{for } |v| < \rho, \\
    f(v) &= a \quad \text{for } \sqrt{5}-\rho < |v| < \sqrt{5} + \rho, \\
    f(v) &= 0 \quad \text{for } |v| > 5, \\
    f(v) &= 1 \quad \text{otherwise}.
\end{align*}

Note that the points $(\pm 2, \pm 1)$ and $(\pm 1, \pm 2)$ belong to the ring where $f = a$. The function $f$ equals $ca^2$ in a small ball around the origin, it equals $1$ in the rest of the ball $B_{\sqrt{5}}$, and zero outside $B_{\sqrt{5}}$.

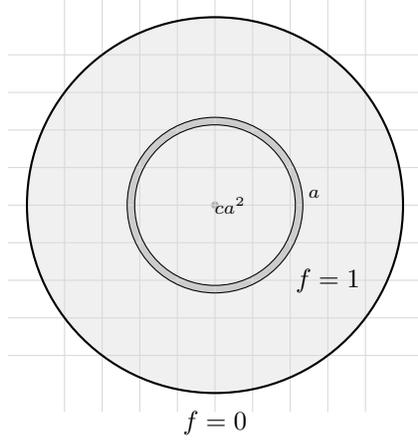
\begin{figure}[h!t]
    \centering
    \begin{tikzpicture}[scale=0.5]
        \foreach \i in {-4,-3,-2,-1,1,2,3,4} {
            \draw[gray!30, thin] (\i,-5.5) -- (\i,5.5);
            \draw[gray!30, thin] (-5.5,\i) -- (5.5,\i);
        }
        \draw[gray!30, thin] (0,-5.5) -- (0,5.5);
        \draw[gray!30, thin] (-5.5,0) -- (5.5,0);

        \fill[gray!12] (0,0) circle (5);

        \fill[gray!40] (0,0) circle ({sqrt(5)+0.1});
        \fill[gray!12] (0,0) circle ({sqrt(5)-0.1});

        \fill[gray!65] (0,0) circle (0.1);

        \begin{scope}
            \clip (0,0) circle (5);
            \foreach \i in {-4,-3,-2,-1,1,2,3,4} {
                \draw[gray!30, thin] (\i,-5.5) -- (\i,5.5);
                \draw[gray!30, thin] (-5.5,\i) -- (5.5,\i);
            }
            \draw[gray!30, thin] (0,-5.5) -- (0,5.5);
            \draw[gray!30, thin] (-5.5,0) -- (5.5,0);
        \end{scope}

        \draw[thick] (0,0) circle (5);
        \draw (0,0) circle ({sqrt(5)+0.1});
        \draw (0,0) circle ({sqrt(5)-0.1});

        \node at (0.4, 0) {\scriptsize $ca^2$};
        \node at (3., -2) {$f=1$};
        \node[right] at ({sqrt(5)}, 0.3) {\scriptsize $a$};
        \node at (0, -5.8) {$f = 0$};

    \end{tikzpicture}
    \caption{Level sets of the function $f$. The dark region near the origin is the ball $B_\rho$ where $f = ca^2$. The gray ring at radius $\sqrt{5}$ is where $f = a$. In the rest of $B_{5}$, $f = 1$.} \label{f:levelsets}
\end{figure}

This is the function that, together with the collision kernel described in the previous subsection, will give us $\partial_t D(f) > 0$. We will prove that if $c$ is small enough, then $\partial_t D(f) \approx a^4 \log a$ for $a \gg 1$.

\begin{remark}
We presented a function $f$ in its simplest form for our computations as a piecewise constant function. It is not difficult to approximate this function by a smooth function and obtain the same result. A similar approximation can be done later for the collision kernel. Thus, the result can be recovered with a smooth function $f$ and a smooth kernel $B$.
\end{remark}

\subsection{Computing the derivative of the entropy production for our function}

We proceed to compute $\partial_t D(f)$ for the function $f$ defined above. We will use the expression for $\partial_t D(f)$ given in Lemma \ref{l:derivative}.

We must analyze $Q$ and $\LL$ for the function $f$ defined above. We will see that $Q$ is of order at most $a^2$ and only for some values of $v$. Moreover, $\LL$ is of order $a^2 \log a$ for some values of $v$. This will allow us to compute the leading order term in $\partial_t D(f)$ and show that it is positive for $a$ sufficiently large.

Let us divide the analysis into estimates for $Q_+$, $Q_-$, and $\LL$.

We write $Q = Q_+ - Q_-$, where $Q_+$ and $Q_-$ are the gain and loss terms in the Boltzmann collision operator, respectively. We have
\begin{align*}
    Q_+(f,f)(v) &= \int_{S_1} f' f'_* d\sigma, \\
    Q_-(f,f)(v) &= f(v) \int_{S_1} f_* d\sigma.
\end{align*}

We observe that for every $v \in \R^2$, $Q_+$ and $Q_-$ cannot be larger than $\approx a^2$. Indeed, $Q_+(v) \approx a^2$ when one of the following two scenarios happen:
\begin{itemize}
\item Both $v'$ and $v'_*$ belong to the ring where $f = a$, or
\item One of $v'$, $v'_*$ belongs to the small ball where $f = c a^2$.
\end{itemize}

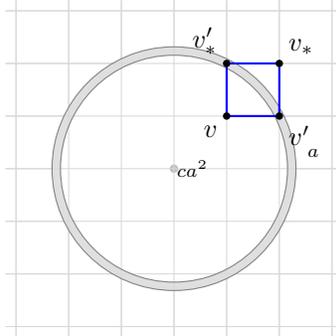
\begin{figure}[ht]
    \centering
    \begin{tikzpicture}[scale=0.7]
        \foreach \i in {-3,-2,-1,1,2,3} {
            \draw[gray!30, thin] (\i,-3.2) -- (\i,3.2);
            \draw[gray!30, thin] (-3.2,\i) -- (3.2,\i);
        }
        \draw[gray!30, thin] (0,-3.2) -- (0,3.2);
        \draw[gray!30, thin] (-3.2,0) -- (3.2,0);

        \fill[gray!25] (0,0) circle ({sqrt(5)+0.08});
        \fill[white!10] (0,0) circle ({sqrt(5)-0.08});

        \fill[gray!50] (0,0) circle (0.08);

        \foreach \i in {-3,-2,-1,0,1,2,3} {
            \draw[gray!30, thin] (\i,-3.2) -- (\i,3.2);
            \draw[gray!30, thin] (-3.2,\i) -- (3.2,\i);
        }

        \draw[gray] (0,0) circle ({sqrt(5)+0.08});
        \draw[gray] (0,0) circle ({sqrt(5)-0.08});

        \draw[thick, blue] (1,1) -- (2,1) -- (2,2) -- (1,2) -- cycle;

        \fill (1,1) circle (2pt);
        \fill (2,1) circle (2pt);
        \fill (2,2) circle (2pt);
        \fill (1,2) circle (2pt);

        \node[below left] at (1,1) {$v$};
        \node[below right] at (2,1) {$v'$};
        \node[above right] at (2,2) {$v_*$};
        \node[above left] at (1,2) {$v'_*$};

        \node at (0.35,0) {\scriptsize $ca^2$};
        \node[right] at ({sqrt(5)+0.12}, 0.3) {\scriptsize $a$};
    \end{tikzpicture}
    \caption{The first scenario for $Q_+$ happens when $f(v') = f(v'_*) = a$. Here $v' = (2,1)$ and $v'_* = (1,2)$. The four points form a square of side length one, with $|v| = \sqrt 2$ and $|v_*| = 2\sqrt 2$.}
\end{figure}

The first scenario takes place for those $v$ that are at distance approximately $\sqrt 2$ or $2\sqrt 2$ from zero. The second scenario takes place for those $v$ that are at distance approximately $1$ from zero. For any other value of $v$, we have $Q_+ \ll a^2$ for large $a$.

The loss term $Q_-$ will be of order $a^2$ when one of the following three scenarios happen:
\begin{itemize}
\item $v$ belongs to the small ball where $f = c a^2$
\item $v_\ast$ belongs to the small ball where $f = c a^2$ somewhere on the domain of integration.
\item Both $v$ and $v_\ast$ belong to the ring where $f = a$.
\end{itemize}

Analyzing these three scenarios, we see that $Q_- \approx a^2$ when either $|v| < \rho$, $|v| \approx \sqrt 2$, or $|v| \approx \sqrt{5}$. For any other value of $v$, we have $Q_- \ll a^2$ for large $a$.

\begin{figure}[ht]
\centering
\begin{minipage}[t]{0.48\textwidth}
\centering
\begin{tikzpicture}[scale=0.5]
    \clip (-4,-4) rectangle (4,4);

    \foreach \i in {-4,-3,-2,-1,1,2,3,4} {
        \draw[gray!30, thin] (\i,-4) -- (\i,4);
        \draw[gray!30, thin] (-4,\i) -- (4,\i);
    }
    \draw[gray!30, thin] (0,-4) -- (0,4);
    \draw[gray!30, thin] (-4,0) -- (4,0);

    \fill[blue!20] (0,0) circle ({2*sqrt(2)+0.1});
    \fill[white] (0,0) circle ({2*sqrt(2)-0.1});

    \fill[blue!20] (0,0) circle ({sqrt(2)+0.1});
    \fill[white] (0,0) circle ({sqrt(2)-0.1});

    \fill[blue!20] (0,0) circle (1.1);
    \fill[white] (0,0) circle (0.9);

    \foreach \i in {-4,-3,-2,-1,0,1,2,3,4} {
        \draw[gray!30, thin] (\i,-4) -- (\i,4);
        \draw[gray!30, thin] (-4,\i) -- (4,\i);
    }

    \draw[gray, dashed] (0,0) circle ({sqrt(5)});

    \draw (0,0) circle (1);
    \draw (0,0) circle ({sqrt(2)});
    \draw (0,0) circle ({2*sqrt(2)});
\end{tikzpicture}
\caption*{$Q_+(v) \approx a^2$}
\end{minipage}%
\hfill
\begin{minipage}[t]{0.48\textwidth}
\centering
\begin{tikzpicture}[scale=0.5]
    \clip (-4,-4) rectangle (4,4);

    \foreach \i in {-4,-3,-2,-1,1,2,3,4} {
        \draw[gray!30, thin] (\i,-4) -- (\i,4);
        \draw[gray!30, thin] (-4,\i) -- (4,\i);
    }
    \draw[gray!30, thin] (0,-4) -- (0,4);
    \draw[gray!30, thin] (-4,0) -- (4,0);

    \fill[blue!20] (0,0) circle ({sqrt(5)+0.1});
    \fill[white] (0,0) circle ({sqrt(5)-0.1});

    \fill[blue!20] (0,0) circle ({sqrt(2)+0.1});
    \fill[white] (0,0) circle ({sqrt(2)-0.1});

    \fill[blue!20] (0,0) circle (0.1);

    \foreach \i in {-4,-3,-2,-1,0,1,2,3,4} {
        \draw[gray!30, thin] (\i,-4) -- (\i,4);
        \draw[gray!30, thin] (-4,\i) -- (4,\i);
    }

    \draw[gray, dashed] (0,0) circle ({sqrt(5)});

    \draw (0,0) circle (0.1);
    \draw (0,0) circle ({sqrt(2)});
\end{tikzpicture}
\caption*{$Q_-(v) \approx a^2$}
\end{minipage}
\caption{The values of $v$ for which $Q_+(v) \approx a^2$ (left) and $Q_-(v) \approx a^2$ (right).}
\end{figure}

Note that the sets where $Q_+$ and $Q_-$ are of order $a^2$ overlap in the ring of radius $\sqrt 2$ with width approximately $\rho$. The value of $Q_-$ there is proportional to the constant $c$, whereas the value of $Q_+$ is unaffected by $c$ (these values also depend on $\rho$, which is fixed at this moment). We pick the constant $c$ so that $Q = Q_+ - Q_- \approx a^2$ is positive in this ring. It is worth noting that the radius $\sqrt 5$ for the ring where $f=a$ was chosen specifically to make these two rings where $Q_+$ and $Q_-$ are of order $a^2$ overlap. This way, we get $Q_+ - Q_- \approx a^2$ where at the same place $Q_- \approx a^2$ and $\LL \approx a^2 \log a$.

Finally, we analyze $\LL$. We will see that $\LL$ is of order $a^2 \log a$ for some values of $v$. This can only happen when $v_\ast$ is able to lie inside $B_\rho$. Thus, it can only happen when $|v| \approx \sqrt 2$. It is the same ring where both $Q_+$ and $Q_-$ are of order $a^2$. In this ring, we have $\LL \approx a^2 \log a$.

From our analysis of $Q$ and $\LL$, we conclude that $Q^2/f$ is of order at most $a^4$ at all points $v \in \R^2$, while $Q \LL$ is of order $a^4 \log a$ for those values of $v$ in the ring of radius $\sqrt 2$ and width approximately $\rho$. As $a \to \infty$, the leading order term in $\partial_t D(f)$, whose expression is given in Lemma \ref{l:derivative}, consists of the integral of $Q \LL$ over the ring where $Q \LL \approx a^4 \log a$. This leading order term is positive, and it dominates the negative term given by the integral of $Q^2/f$, which is of order at most $a^4$. Therefore, $\partial_t D(f) > 0$ for $a$ sufficiently large.

This concludes the proof of Theorem \ref{t:main}.

\bibliographystyle{plain}
\bibliography{ed-increasing}

\end{document}